\documentclass[11pt]{article}

\usepackage{graphicx}
\usepackage{latexsym}

\usepackage[letterpaper,width=135mm,height=203mm]{geometry}

\newtheorem{theorem}{Theorem}
\newtheorem{lemma}[theorem]{Lemma}
\newtheorem{conjecture}[theorem]{Conjecture}
\newtheorem{proposition}[theorem]{Proposition}
\newtheorem{observation}[theorem]{Observation}

\newenvironment{Proof}{\proofing}{\QED}
\newcommand{\QED}{\hspace{8mm}\mbox{\textsc{qed}}\smallskip}
\newcommand{\proofing}{\textsc{Proof.}}

\newcommand{\dd}[1]{\textbf{\textit{#1}}}

\newcommand{\curlyP}{\mathcal{P}}
\newcommand{\Cor}{\mathop{\mathit{Cor}}}
\newcommand{\diam}{\mathop{\mathrm{diam}}}
\newcommand{\Ttilde}{\tilde{T}}
\newcommand{\Itilde}{\tilde{I}}
\newcommand{\Ptilde}{\tilde{P}}

\makeatletter
\long\def\@caption#1[#2]#3{\begingroup \@parboxrestore 
\if@minipage \@setminipage \fi \normalsize \sffamily \@makecaption {\csname fnum@#1\endcsname }{\ignorespaces #3}\par \endgroup}
\makeatother

\begin{document}

\nocite{*}

\begin{center}
{\Large \textbf{On the Relationships between \\[1mm] 
Domination, Isolation, and Packing}} \\[4mm]
Geoffrey Boyer\textsuperscript{1}, 
Wayne Goddard\textsuperscript{2,3},
Michael A. Henning\textsuperscript{3} \\[2mm]
\textsuperscript{1}Department of Mathematics and Statistics, College of Charleston \\
\textsuperscript{2}School of Mathematical and Statistical Sciences, Clemson University \\
\textsuperscript{3}Department of Mathematics and Applied Mathematics, University of Johannesburg
\end{center}

\begin{abstract}
We consider the relationships between the domination number of graph, denoted $\gamma$, and the distance-$2$ domination number, denoted $\gamma_2$, and three parameters that lie between them: the packing number, denoted $\rho$, the lower packing number, denoted $\rho_L$,
and the isolation number, denoted $\iota$. There has been recent attention on the question of whether $\gamma/\rho$ is bounded or unbounded for various families of graphs. We consider similar questions for the ratios of the five parameters. In particular we show that, while $\gamma/\rho_L$ is unbounded in trees, it holds that $\iota/\gamma_2$ is less than $2$ for all trees. Further, $\gamma/\rho_L$ is
at most $3$ in interval graphs, at most~$4$ in permutation graphs, and at most $5$ in general asteroidal-triple-free graphs. We also show that every tree has a set of vertices that is both isolating and a packing,
and characterize trees where $\rho=\rho_L$.
\end{abstract}


\section{Introduction}

Recently there has been renewed interest in the relationship between the domination number and packing number (also called $2$-packing number or $2$-independence number)
of a graph. In this paper we contribute several more results, as well as consider the relationships involving the isolation number and packing number of a graph.
All our graphs are finite, simple, and undirected.

There is a long history on the topic. 
Perhaps the earliest connection was provided by Meir and Moon~\cite{MM75tree} who
showed that trees have equal packing and domination numbers. Later this equality was extended to strongly chordal graphs by Farber~\cite{Farber84chordal}.
In~\cite{HLR11pack} the current avenue was triggered by the observation that the domination
number is at most the product of the packing number and the maximum degree.
Continuing down the road of degree-bounded graphs, for subcubic claw-free graphs they showed
that the domination number is at most twice the packing number. 
Another recent direction is replacing domination by independent domination: for example,
Cho and Kim~\cite{CK24idomPack} proved that the independent domination number is at most three times the packing number in subcubic graphs. 

There has been continued work on the ratio of the domination and packing numbers.
G\'omez and Guti\'errez \cite{GG25domPack} focussed on bipartite cubic graphs and showed that the ratio of the 
domination number to the packing number in such a graph is at most $120/49$.
Bonamy et al.~\cite{BCGY25} investigated graph classes where the ratio is bounded by 
an absolute constant. Inter alia, they showed that 
$2$-degenerate, asteroidal-triple-free, and unit-disk graphs are classes with bounded ratio.
They also improved the known ratio for planar graphs and graphs of bounded treewidth,
while providing further instances of graph classes where the ratio is unbounded.
Some unresolved cases were resolved by D\'ucz and Gujgiczer~\cite{DG26}.

A less well-known parameter is the isolation number of a graph.
Isolation was defined and studied by Lewis et al.~\cite{LHHF10vertexEdge}
and by Caro and Hansberg~\cite{CH}. More results and history are to be found in \cite{BG24first}. 
Much of the work on the isolation number has focussed on bounds for classes. For example, Canales et al.~\cite{CHMM15distance} and 
Tokunaga et al.~\cite{TJK19outerplanar} considered maximal outerplanar graphs, 
Ziemann and \.{Z}yli\'{n}ski~\cite{ZZ20cubic} considered cubic graphs, and 
Adhya et al.~\cite{AMB26perm} considered permutation graphs.
Caro and Hansberg~\cite{CH} also provided a generalization of isolation, which has been the subject
of multiple recent papers, but we do not explore it here.

We proceed as follows. 
In Section~\ref{s:prelim} we provide the necessary definitions and notation.
In Section~\ref{s:packIso} we study graphs that have a set that is both an isolating set and a packing, showing that every tree has such a set.
In Section~\ref{s:lowerPack} we consider two questions about maximal packings in trees. 
In Section~\ref{s:regular} we provide some observations about how these parameters relate in cubic and degree-bounded graphs.
In Section~\ref{s:perfect} we extend some of the earlier work to incorporate the isolation and lower packing numbers.
We conclude in Section~\ref{s:conclude} with some further thoughts.


\section{Preliminaries and Definitions} \label{s:prelim}

The \dd{closed neighborhood} of a set $S$ of vertices, which we denote by $N[S]$, is the set consisting of $S$ and all vertices with a neighbor in $S$. 
A \dd{dominating set} of a graph $G$ is a set $S$ of vertices such that every vertex of $G$ is in $N[S]$; 
the \dd{domination number} is the minimum size of a dominating set and is denoted $\gamma(G)$. 
An \dd{isolating set} of $G$ is a set $S$ of vertices such that $G-N[S]$ has no edge; 
the \dd{isolation number} is the  minimum size of an isolating set and is denoted~$\iota(G)$. 
This parameter is also called the \dd{vertex-edge} domination number  
(where it is denoted $\gamma_{ve}(G)$). 
A \dd{packing}, also called a \dd{$2$-packing}, is a set $S$ of vertices such 
that the closed neighborhoods of each vertex in $S$ are disjoint. The \dd{packing number} is the 
maximum size of a packing and is denoted $\rho(G)$. We will also consider the 
\dd{lower packing number}, denoted $\rho_L(G)$, which is the minimum size of a maximal packing.
Furthermore, a set $S$ is \dd{distance-$2$ dominating} if every vertex of $G$ is within distance $2$ of some
vertex of $S$. The \dd{distance-$2$ domination number}, denoted here by $\gamma_2(G)$, is the minimum 
size of a distance-$2$ dominating set. 

It is well known and immediate that the domination number is at least the packing number, since a dominating set must contain at least one vertex from the closed
neighborhood of each vertex in a packing. By definition, the domination number is at least the isolation number. Furthermore, any isolating set in a
connected nontrivial graph is distance-$2$ dominating, as is every maximal packing. Thus we have for a connected nontrivial graph~$G$:
\[
      \gamma_2(G) \le \left\{ \begin{array}{cc} \rho_L(G) \le \rho(G) \\ \iota(G) \end{array} \right\} \le \gamma(G) .
\]                          

In general the isolation number
and the packing number are incomparable. For example, the rooks graph $K_m \Box K_m$ has packing number (and lower packing number)  $1$ and isolation number $m-1$. 
In contrast, the octopus $O_m$ formed by subdividing every edge of $K_{1,m}$ has isolation number $1$ and packing number $m$.
The octopus has lower packing number $1$. But consider 
the tree $D_m$ obtained by taking two copies of $O_m$ and adding an edge joining their centers;
this has isolation number $2$ and lower packing number $m$. 

We will need the following result of Lewis et al. and Caro and Hansberg.

\begin{proposition}  \label{p:isoLower} \cite{LHHF10vertexEdge,CH}
For any graph $G$ it holds that  $\iota(G) \ge m / \Delta^2$, where $m$ is the number of edges and $\Delta$ the maximum degree.
\end{proposition}

\section{Graphs with a Packing Isolating Set} \label{s:packIso}

Given that the isolation number and packing number are in general incomparable, one topic of interest is graphs 
that have a set that is both a packing and an isolating set (what we call a \dd{packing isolating set}). 
Graphs that have dominating sets that are packings have been considered under several names, including \dd{efficient domination} and \dd{perfect domination}. 
Most graphs do not have such sets, and even among trees there are examples with no such set (the double-star is an example). 
As regards isolation, most graphs do not have a packing isolating set either ($C_5$ is an example). 
But we show here that every tree has a packing isolating set. 
Note that questions about sets that are both a $2$-dominating set and a packing were recently explored by Bujtás at al.~\cite{BCKZ26}.

By the definitions, if a graph $G$ has  a  packing isolating set $S$ then $\iota(G) \leq |S| \leq \rho(G)$. 
Further, such a set $S$ is necessarily a minimal isolating set. 
By a result in \cite{LHHF10vertexEdge} it follows that $|S|\leq n/2$, where $n$ is the order of $G$. It is possible for 
there to be equality. 
Define the \dd{corona} $\Cor(G)$ of a graph $G$ as the graph obtained by adding one leaf vertex adjacent to each vertex. 
It is immediate that the set of leaf vertices in a corona is a packing isolating set.  

\subsection{Trees}

Apart from the work of Meir and Moon~\cite{MM75tree} there has been other work on packings in trees;
bounds on various versions of packing were considered in~\cite{henningPack} and more recent work on 
``open packings'' by Hartnell and Rall~\cite{HR20open}. We consider here the relationship between packings and isolating sets
in trees. We observe first that is not true, even in trees, that a maximum packing is guaranteed to be an isolating set. Consider for example
the set consisting of the two ends of the path~$P_6$ on six vertices. But we show here that 
every tree has a packing isolating set.

\newpage

\begin{theorem}
Every tree $T$ has a set that is both a packing and an isolating set.
\end{theorem}
\begin{Proof}
The proof by induction on the order of $T$. 
Trivially a central vertex of a tree of diameter at most $3$ forms a packing isolating set. So we may assume the diameter is at least~$4$.
Consider a longest path in the tree, say ending $wxyz$. Let $M = N[y]-\{x\}$; 
note that all vertices of $M$ except $y$ are leaf vertices of $T$.
Consider the tree $T' = T - M$. By the inductive hypothesis the tree $T'$ has a packing isolating
set $J$. 

If $x \in J$, then $J$ is an isolating set of $T$, and we are done. 
So assume $x \notin J$. If $J$ dominates $x$, then adding $z$ to $J$ produces
an isolating set of $T$, and it remains a packing.
On the other hand, if $x$ is not dominated by $J$, then adding $y$ to~$J$ produces
an isolating set of $T$, and it remains a packing.
The result follows. 
\end{Proof}

\begin{figure}[h]
\noindent\hfill
\includegraphics[scale=0.9]{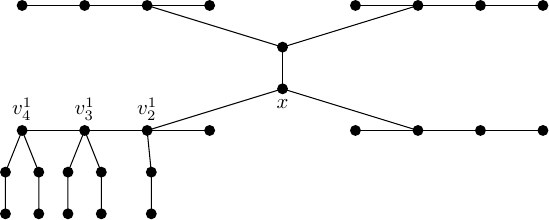}
\hfill\null
\caption{A tree $\tilde{T}$ with packing number more than isolation number.}
\label{f:treePack}
\end{figure}

Nevertheless, there exist trees where no minimum isolating set is a packing 
and no maximum packing is isolating. 
Consider for example the following tree $\Ttilde$ (inspired by an example in~\cite{BCM25well}).
Take four disjoint copies of $P_4$;  say the $i$\textsuperscript{th} copy is  $v^i_1,v^i_2,v^i_3,v^i_4$.
Then add a $K_2$ with vertices $x$ and $y$ and join $x$ to $v^1_2$ and $v^2_2$ and join $y$ to $v^3_2$ and $v^4_2$.
Finally add two copies of~$K_2$ adjacent to $v^1_3$, two copies of $K_2$ adjacent to $v^1_4$, and one copy of $K_2$ adjacent
to $v^1_2$. See Figure~\ref{f:treePack}.
It is easy to see that every minimum isolating set of the result contains both  $v^1_3$ and $v^1_4$, and is thus not a packing;
on the other hand, every maximum packing contains each $v^i_1$ and thus does not dominate the edge $xy$.
In particular, it follows that for every packing isolating set $S$ it holds that $\iota(\Ttilde) < |S| < \rho(\Ttilde)$.

While we observed earlier that a packing isolating set (of a graph with minimum degree at least $1$) contains at most half the vertices, 
trees do not necessarily have two disjoint packing isolating sets. For example, 
let tree $F$ be formed by subdividing each edge of $K_{1,3}$; and form 
tree $T$ by taking three copies of $F$
and add two edges
between the centers of $F$. It suffices to note that while $F$ itself
does have two disjoint packing isolating sets, one of these sets must contain the 
center vertex. 

\subsection{Graphs Without a Packing Isolating Set}

\begin{figure}[h]
\noindent\hfill
\includegraphics[scale=0.9]{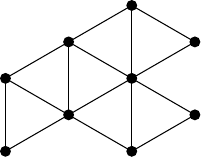} \qquad\qquad\quad
\includegraphics[scale=0.9]{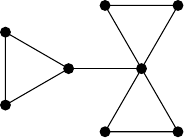}
\hfill\null
\caption{Two graphs with no packing isolating set.}
\label{f:noPackIso}
\end{figure}

The existence of a packing isolating set, shown earlier for the family of trees, does not seem to extend to other families of graphs.
For example, there exist maximal outerplanar graphs (MOPs) that have no packing isolating set---an example is given on the left in Figure~\ref{f:noPackIso}---
and interval graphs that have no packing isolating set---an example is given on the right in Figure~\ref{f:noPackIso}.
Furthermore, it is not true that every bipartite graph has a packing isolating set. 
Take for example the Mobius ladder of order 10 (the circulant $C_{10}[1,5]$). 

Indeed:

\begin{theorem}
Testing whether a graph has a packing isolating set is NP-hard.
\end{theorem}
\begin{Proof}
We reduce from 1-in-3 SAT. 
Given a boolean formula $\phi$, we build a graph $G_\phi$ such that 
there exists a packing isolating set of $G_\phi$ if and only if there is a $1$-in-$3$-satisfying assignment of $\phi$.
We construct $G_\phi$ as follows (similar to the original proof that Domination is NP-hard). 
For each variable $v$ in $\phi$, we introduce a $K_2$ with vertices labeled $v$ and not-$v$.
For each clause $c$ in $\phi$, we introduce a $C_5$ with one designated vertex $w_c$. 
Then we join each $w_c$ to the three literals contained in $c$.
 
Suppose there exists a packing isolating set $J$ of $G_\phi$. 
Since the graph $C_5$ has no packing isolating set, the set $J$ cannot contain any of the $w_c$.
Thus $J$ must contain one literal from each variable and $J$ must dominate all of the $w_c$. 
Furthermore, any two literals in $J$ must be at distance more than $2$,
and thus each $w_c$ is dominated by exactly one vertex of $J$. It follows that the set of literals in $J$ represents a 
1-in-3 satisfying assignment for $\phi$.

Conversely, given a satisfying assignment, create set $J$ by taking the true literals and one vertex from each clause $c$
that is not adjacent to $w_c$. This is a packing isolating set of $G_\phi$.
\end{Proof}


\section{More about Maximal Packings in Trees} \label{s:lowerPack}

\subsection{A Bound For Isolation but not Domination}

For trees, Meir and Moon~\cite{MM75tree} showed that the packing and domination numbers are equal.
In contrast, the ratio $\gamma(T) / \rho_L(T)$ is unbounded: consider for example the octopus $O_m$ from earlier.

Turning to isolation, it follows that in trees the isolation number is at most the packing number.
Further, as noted above, the ratio $\iota(T)/\rho(T)$ can be arbitrarily small.
But we observe here that, unlike in general graphs, 
in trees the ratio $\iota/ \rho_L$ cannot be arbitrarily large:

\begin{theorem}  \label{t:iotaRhoL}
For any tree $T$ it holds that  $\iota(T) \le 2 \rho_L(T) -1 $.
\end{theorem}

In fact we show that $\iota(T) \le 2 \gamma_2(T) -1 $, and the theorem follows  
since a maximal packing is a distance-$2$ dominating set. 

\begin{lemma} \label{l:treeeExtends}
For any tree $T$ and any distance-$2$ dominating set $P$ of $T$, the set $P$ can be extended to 
an isolating set of $T$ by adding at most $|P|-1$ vertices.
\end{lemma}
\begin{Proof}
The proof is by induction on the order of $T$. Assume first there is a leaf vertex $\ell$ 
of $T$ that is not in $P$. Then vertex $\ell$ must be at
distance at most $2$ from $P$ and so its neighbor is dominated by $P$.
Let $T_0$ be the tree formed from $T$ by deleting $\ell$. By the inductive hypothesis we can 
extend $P$ to an isolating set of $T_0$ by adding at most $|P|-1$ vertices; by the
discussion this set is also an isolating set of $T$. Thus we may assume all leaf vertices
of $T$ are in $P$.

Assume second there is a vertex $p$ of $P$ that is not a leaf vertex of $T$. Then 
there exist two subtrees $T_1$ and $T_2$ of $T$ such that every vertex of $T$ is in at 
least one of the subtrees, and only $p$ is in both. Let $P_i$ be the restriction
of $P$ to $T_i$. By the inductive hypothesis we can extend $P_i$ to an isolating
set of $T_i$ by adding at most $|P_i|-1$ vertices. The union of these two sets
is an isolating set of $T$, and it contains $P$ and at most $|P_1|-1+|P_2|-1 = |P|-1$ other vertices, as required.
Thus we may assume all vertices of $P$ are leaf vertices of $T$.

Root the tree $T$ at some leaf vertex $\ell$. If $T$ is a path, then by distance domination of~$P$
the tree $T$ has at most six vertices; the set $P$ can be extended to an isolating set by adding
one of the central vertices. So we may assume $T$ is not a path. Let $x$ be a vertex
of degree more than $2$ farthest from $\ell$. Let $z$ be a descendant of $x$ 
farthest from $x$. Let $y$ be the child of $x$ on the path to $z$; possible $y=z$.
By the choice of $x$, 
vertex $y$ and its descendants all have degree at most $2$ in~$T$. Define $T'$ as
the tree obtained from $T$ by deleting $y$ and all its descendants, and define
$P' = P - \{ z \}$. The set $P'$ is a distance-$2$ dominating set of $T'$: the remaining 
leaf descendants of $x$ are no farther than $z$, so the removal of $z$ cannot break the 
distance domination. By the inductive hypothesis the set $P'$ can be extended to 
an isolating set of $T'$ by adding at most $|P'|-1$ vertices. This isolating
set can be extended to an isolating set of $T$ by adding $y$ and $z$, for a total of $|P|-1$ vertices other than $P$, as required.
\end{Proof}

\begin{figure}[h]
\noindent\hfill
\includegraphics{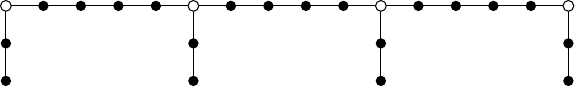}
\hfill\null
\caption{The tree $P_4^*$ with $\iota=7$ and $\gamma_2  = \rho_L=4$.}
\label{f:iotaRhoL}
\end{figure}

As a consequence of Lemma~\ref{l:treeeExtends}
we obtain  Theorem~\ref{t:iotaRhoL}.
There are trees where the bound of Theorem~\ref{t:iotaRhoL} is attained. Consider the following construction.
Starting with any graph~$H$, form the graph~$H^*$ by first creating the corona $\Cor(H)$, 
then subdividing each original edge with four vertices, and subdividing each corona edge with one vertex. For example,
$P_4^*$ is drawn in Figure~\ref{f:iotaRhoL}.

\begin{lemma}
For any graph $H$ with $n$ vertices and $m$ edges it holds that $\rho_L( H^* ) = n$ and $\iota( H^*) = n+m$.
\end{lemma}
\begin{Proof}
The vertices of $H$ form a maximal packing of $H^*$, and every maximal packing of $H$ contains a vertex within
distance two of each vertex of $H$. Thus $\rho_L(H^*) = n$.
On the other hand, to form an isolating set of $H^*$ we need to dominate each leaf-edge; so without loss of generality the isolating set 
contains every vertex of $H$. But an additional vertex is needed to dominate the edge midway between each adjacent pair
of vertices of $H$. Hence $\iota( H^* ) = n+m$.
\end{Proof}

\noindent
Thus for any tree $H$ the tree $H^*$ is an example of equality in Theorem~\ref{t:iotaRhoL}.

\begin{figure}[h]
\noindent\hfill
\includegraphics{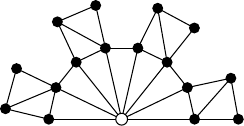}
\hfill\null
\caption{A MOP with $\iota \gg \rho_L$}
\label{f:mop}
\end{figure}

The boundedness of $\iota/\rho_L$ does not generalize to other graph families. For example, 
if $H$ is any bipartite graph then $H^*$ is bipartite and by choosing $H$ to be the complete bipartite graph
one obtains a graph with $\iota/\rho_L$ arbitrarily large.
Further, one can build a $2$-tree, indeed a MOP, with a large ratio.
For $m\ge 1$ build MOP $M_m$ as follows.
Start with the path $v_1 \ldots v_{2m}$. 
Add one vertex adjacent to all vertices on the path.
Further, for $1\le j \le m$ add two new adjacent vertices $x_i$ and $y_i$ such that $x_i$ is adjacent to 
$v_{2i}$ and $y_i$ is adjacent to both $v_{2i}$ and $v_{2i+1}$. The MOP $M_4$ is drawn in Figure~\ref{f:mop}.
It can be shown that $\rho_L( M_m) =1$ but
$\iota(M_m) = m$.

\subsection{Well-Packed Trees}

One question of interest with any collection of sets in a graph is when do all the sets have the same size. Thus, 
for example, in well-covered graphs  all maximal independent sets have the same size and in well-dominated graphs
all minimal dominating sets have the same size.
Recently, Boutrig, Chellali, and Meddah~\cite{BCM25well} determined the trees where every minimal independent isolating set has the same size
(what they called well-$ve$-covered). 
Also, Büyükçolaka~\cite{Büyükçolaka} provided a characterization of the trees where every minimal isolating set has the same size (what they called well-$ve$-dominated).
In fact the families coincide. As Büyükçolaka observed, the condition of being well-$ve$-dominated is a stronger condition than being well-$ve$-covered,
and thus the set of well-$ve$-dominated graphs is necessarily a subset of the well $ve$-covered graphs.
But it easy to check that all the trees determined in~\cite{BCM25well} are also well-$ve$-dominated.

\begin{figure}[h]
\noindent\hfill
\includegraphics{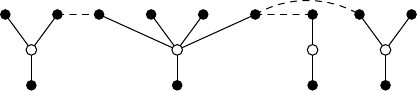}
\hfill\null
\caption{A tree in $\curlyP$}
\label{f:wellPack}
\end{figure}

We consider here the question for packing number. (We note that the question for ``open packing number'' was resolved by Hartnell and Rall~\cite{HR20open}.)
We say a graph $G$ is \dd{well-packed} if every maximal packing has the same size; equivalently $\rho_L(G) = \rho(G)$.
Define $\curlyP$ as the set of all trees constructed as follows. Start with a disjoint collection of $s$ stars each of order at least $3$; let $Y$ 
be the set of all the centers. Then add $s-1$ edges to make the collection into a tree such that no added edge is incident to $Y$ and also every vertex of 
$Y$ still has at least one leaf neighbor. Equivalently, consider any tree where the support vertices (meaning those with leaf neighbors) 
form a set $Y$ that is both a packing and a dominating set. An example is shown in Figure~\ref{f:wellPack}.

We note that  the well-$ve$-dominated / well-$ve$-covered trees are the trees in~$\curlyP$ where
every vertex of $X$ has at most one non-leaf neighbor.

\begin{theorem}   
For $n\ge 3$, a tree $T$ is well-packed if and only if $T$ is in $\curlyP$.
\end{theorem}
\begin{Proof}
(1) We need to show that every tree in $\curlyP$ is well-packed. Consider a tree $T \in \curlyP$.
Every packing of $T$ contains at most one vertex from each constituent star.
Furthermore, the leaf vertex is at least distance $3$ from every other star, and so can be added to any packing that contains no
vertex from that star. In other words, every maximal packing of $T$ contains exactly one vertex from each constituent star.
This means that $T$ is well-packed.\smallskip

(2) Assume $T$ is well-packed. We prove that $T \in \curlyP$ by induction on its order.
If $T$ is a star then it is in $\curlyP$. So assume its diameter is at least $3$.
Consider a longest path in $T$, say ending $wxyz$.
By the maximality of the path, every neighbor of $x$ except (possibly) $w$ is either a leaf or 
all of its neighbors except $x$ are leaves.
Consider a maximal packing $P_x$ of $T$ that contains vertex $x$. 
If vertex $x$ has degree more than $2$ in $T$, then one can replace $x$ in $P_x$ by $z$ and one neighbor other than $w$ or $y$,
and have a larger packing, a contradiction of $T$ being well-packed. It follows that $x$ has degree~$2$ in $T$.

Let $T'$ be the subtree of $T$ consisting of all vertices that are closer to $w$ than to~$x$. The tree $T'$ must be
well-packed, since each maximal packing of $T'$ can be extended to one of $T$ by 
adding vertex $z$. 
If $T'$ has order $1$ or $2$, then both $\{x\}$ and $\{z,w\}$ are maximal packings of $T$, a contradiction of $T$ being well-packed.
So it must be that $T'$ has order at least $3$.
By the inductive hypothesis, it follows that the tree $T'$ is in~$\curlyP$. 

Let $Y'$ be the centers of the constituent stars of $T'$. To show that $T$ is in $\curlyP$ it suffices to show that in $T$ the set
$Y' \cup \{ y \}$ is a packing and every vertex thereof has a leaf neighbor.

Suppose $w$ is in $Y'$; say with a leaf neighbor $u$.
A maximal packing of $T'$ containing $u$ can be extended to one of $T$ by adding $z$, or by removing
$u$ and adding $x$. This yields maximal packings of $T$ of two different sizes, a contradiction of $T$ being well-packed. 
It follows that $w\notin Y'$. In particular, $Y' \cup \{ y \}$ is a maximal packing of $T$.
By the construction of $\curlyP$ there is a unique neighbor of $w$ in $T'$ that is in $Y'$, say $v$.

Suppose $w$ is the \textbf{only} leaf-neighbor of vertex $v$ in~$T'$ (meaning $v$ does not have a leaf neighbor in $T$). 
Then in $T$ build the set $B$ by adding, for each neighbor of $v$ 
a neighbor of that vertex that is not $v$. (Thus for example $B$ contains $x$.)
The set $B$ forms a packing; hence it can be extended to a maximal packing $B^*$ 
of $T$. In order for $B^*$ to have the same size as $Y' \cup \{ y \}$ it must contain a vertex of $N[v]$. 
But by construction every such vertex is within distance two of a vertex of $B$, a contradiction. 
It follows that $v$ has a leaf neighbor in $T$, whence $T$ is in~$\curlyP$.
\end{Proof}

There does not seem to be a simple description of well-packed graphs in more general families, e.g.\ well-packed MOPs.

\section{Regular and Degree-Bounded Graphs}  \label{s:regular}

We consider next graphs with bounded maximum degree, and especially cubic graphs.

\subsection{Construction}

\begin{figure}[h]
\noindent\hfill
\includegraphics{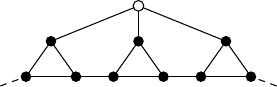} \qquad 
\includegraphics{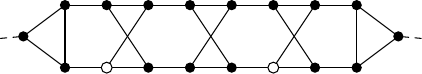}
\hfill\null
\caption{The subcubic graphs $H_{10}$ and $H_{18}$}
\label{f:h1018}
\end{figure}

We will need the following construction, which we call a \dd{bracelet}.
Assume $F$ is a connected graph with two vertices of degree $2$ and all other vertices of degree $3$. 
Then we let $B_s(F)$ denote a connected cubic graph obtained
from $s$ disjoint copies of $F$ and adding $s$ edges. If there is an automorphism of $F$ that interchanges its two degree-$2$ vertices,
then $B_s(F)$ is unique up to isomorphism. 
Let $H_6 = K_{3,3}-e$.
Let graph $H_{10}$ be the subcubic graph drawn on the left in Figure~\ref{f:h1018}.
Let graph $H_{18}$ be the subcubic graph drawn on the right in Figure~\ref{f:h1018}.

\begin{lemma} \label{l:bracelet}
For $s\ge 1$: \\
(a) It holds that $\iota(B_s(H_6))=s$ and $\gamma(B_s(H_6)) = 2s$. \\
(b) It holds that $\rho_L( B_s( H_{10} ) ) =   \gamma_2 ( B_s( H_{10} ) ) = s$ and 
$\gamma( B_s( H_{10} ) ) = 3s$. \\
(c) It holds that $\rho_L(B_s(H_{18} ) ) = 2s $ and $\iota( B_s(H_{18} ) ) = 5s$.
\end{lemma}
\begin{Proof}
(a) 
Any set consisting of the same vertex of each copy of $H_6$ is isolating. 
On the other hand, any dominating set must contain at least two vertices from each copy of $H_6$, and it is easy to 
find one with size $2s$.

(b)
It is immediate that the set of hollow vertices from each copy of $H_{10}$ forms a maximal packing 
of $B_s( H_{10} )$, while every distance-$2$ dominating
set must contain at least one vertex from each copy of $H_{10}$. 
On the other hand, it can be checked that every dominating set contains at least
three vertices from each copy of $H_{10}$. 

(c) The set of the two hollow vertices from each copy form a maximal packing. 
And one can show (or check by computer) that $H_{18}$ minus the two degree-two vertices has isolation number $5$, 
and thus five vertices are needed from each copy of $H_{18}$ to make an isolating set.
\end{Proof}

\subsection{Bounds on Ratios}

In~\cite{HLR11pack} it was noted that 
if $G$ is nontrivial connected graph then $\gamma(G) \le \Delta(G) \rho(G)$. In fact the same reasoning immediately provides the following bound
(though we could not find a record of this in the literature).

\begin{observation}  \label{ob:gammaGamma2}
If $G$ is nontrivial connected graph, then $\gamma(G) \le \Delta(G) \gamma_2 (G)$. 
\end{observation}

It is not hard to find examples of equality, and indeed in the implication $\gamma(G) \le \Delta(G) \rho_L (G)$. By 
Lemma~\ref{l:bracelet} the bracelet $B_s( H_{10} )$ is an example.
In contrast, graphs $G$ where $\gamma(G) = \Delta(G) \rho(G)$ seem rare. Indeed it is noted in~\cite{HLR11pack} that such graphs must
be regular and well-packed, and there are exactly three cubic examples. Moreover, they 
conjecture that $\gamma(G) \le 2 \rho(G) + 1$ for all cubic graphs. 
We note that an upper bound on the ratio $\gamma(G)/\rho(G)$ can be deduced by combining a lower bound for $\rho(G)$ and an
upper bound for $\gamma(G)$. Consider for example connected cubic graphs of order $n$. 
Kostochka and Stocker~\cite{KS09cubic} showed that for $n>8$ that $\gamma(G) \le \frac{5}{14}n$; we~\cite{GHpack} showed
that $\rho(G) \ge \frac{17}{132} (n - 3)$. This implies that for all $n$ sufficiently large the ratio $\gamma/\rho$ is at most $2.78$.

Turning to isolation, Observation~\ref{ob:gammaGamma2} also implies for example that $\iota(G) \le \Delta(G) \rho_L(G)$.
But here it seems that equality is not achievable for cubic graphs other than the Petersen graph.
Except for that graph, the largest example ratio $\iota(G)/ \rho_L(G)$ in cubic graphs we know is $\frac{5}{2}$, achieved
by $B_s(H_{18})$ (see Lemma~\ref{l:bracelet}).

Observation~\ref{ob:gammaGamma2} also implies for example that $\gamma(G) \le \Delta(G) \iota (G)$. 
It is not hard to find examples of equality in this bound. For example, for 
subcubic graphs,
one construction is to take the corona $\Cor( C_m)$ where $m$ is a multiple of $3$.
Then $\gamma=m$ but the set consisting of every third vertex on the cycle shows that $\iota = m/3$.
Or more generally, take the corona $\Cor(H) $ of an $r$-regular graph $H$ of order $m$ that has an efficient dominating set (that is,
a dominating set that is a packing). As observed for example in \cite{BMMS},
the isolation number of the corona is the domination number of the base. Thus $\gamma( \Cor(H) ) = m = (r+1) \gamma(H) = (r+1) \iota ( \Cor(H) )$.

Nevertheless, for regular graphs the ratio $\gamma(G) / \iota(G)$ is definitely smaller, at least for cubic graphs.
By Proposition~\ref{p:isoLower} it holds that $\iota(G) \ge n/6$ for a cubic graph.
By the result of Kostochka and Stocker~\cite{KS09cubic},  a cubic graph has domination number at most $5n/14$ (for $n>8$).
The bounds on the isolation and domination numbers yield that in cubic graph the ratio $\gamma(G) / \iota(G)$  is at most~$\frac{15}{7}$.
On the other hand, we believe:

\begin{conjecture}
For a cubic graph $G$ it holds that $\gamma(G) \le 2 \iota(G)$.
\end{conjecture}

\noindent
There are example cubic graphs $G$ where $\gamma(G) = 2 \iota(G)$. Namely, $B_s( H_6 )$.

A weaker version of the upper bound in Observation~\ref{ob:gammaGamma2} is that $\iota(G) \le \Delta(G) \rho(G)$. 
For maximum degree $3$, and indeed cubic graphs, there is an example of equality, namely the Petersen graph. 
But we wonder whether in connected subcubic graphs, except the Petersen graph, that $\iota(G) \le 2 \rho (G)$.
(There are multiple small subcubic graphs where $\iota/\rho=2$.)


\section{Families of Perfect Graphs} \label{s:perfect}

Finally we consider bounds on ratios of parameters in some families of perfect graphs.
For any established bound $\gamma\leq c\rho$, it is immediate that $\iota\leq c\rho$ since $\iota\leq\gamma$. 
However, if $\gamma/\rho$ is unbounded, this does not imply the same is true for $\iota/\rho$. We note that there 
are cases where this implication is almost immediate, but also cases where in fact $\iota/\rho$ is bounded.

As regards bipartite graphs, Bonamy et al.~\cite{BCGY25} point out that the ratio $\gamma(G) / \rho (G)$ is unbounded. It is to be noted that 
the actual example they provide is mistakenly not bipartite. But one can readily provide such an example using ideas there that 
extends to isolation. Consider the incidence graph $I_k$ of a projective plane of order $k$. 
Since $I_k$ is  bipartite with diameter $3$, it has packing number $2$. But its domination number is $2k$ (see \cite{GRMdesign}).
Indeed by Proposition~\ref{p:isoLower}, the isolation number of $I_k$ is at least $(k^2+k+1)/(k+1)$,
and thus $\iota/\rho$ can be arbitrarily large in bipartite graphs. 

As regards split graphs, Bonamy et al.~\cite{BCGY25} point out that the ratio $\gamma(G) / \rho (G)$ is unbounded for split graphs and hence for chordal graphs.
However, for isolation, all split graphs have isolation number $1$. We note here that it does remain the case that the ratio $\iota/\rho$ is unbounded in chordal graphs.
We use a special case of the lexicographic product (a more general result is given in~\cite{BDJKRTproduct}). The following is straight-forward:

\begin{observation}
If $G[K_2]$ denotes the graph obtained from graph $G$ by duplicating each vertex to form a $K_2$, then $\iota(G[K_2]) = \gamma(G[K_2]) = \gamma(G)$.
\end{observation}

Since $G[K_2]$ remains chordal if $G$ is, and $\rho(G[K_2]) = \rho(G)$, the unboundedness of $\iota/\rho$ in chordal graphs follows.

\begin{figure}[h]
\noindent\hfill
\includegraphics{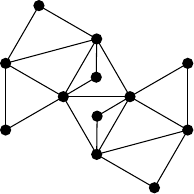}
\hfill\null
\caption{A $2$-tree with $\gamma=2\rho$}
\label{f:2tree}
\end{figure}

As regards $k$-trees, 
Bonamy et al.~\cite{BCGY25} showed that $\gamma(G) / \rho(G) \le k$ for a graph with treewidth $k$, improving and extending
the bound of $3$ provided for MOPs by~\cite{GG25domPack}. In the latter paper it is observed that the sun graph of order $6$ is
an example where $\gamma/\rho = 2$, and they conjecture that this is the only such MOP. As regards $2$-trees,
there are more examples where the ratio $\gamma/\rho$ is $2$. For example, an infinite family is obtained by taking the sun graph and 
duplicating each of the degree-$2$ vertices independently any number of times. 
But there are other examples that do not have twins. Figure~\ref{f:2tree} shows the one of smallest order.

As regards $k$-degenerate graphs,
Bonamy et al.~\cite{BCGY25} showed that $\gamma(G) / \rho(G)$ is unbounded for $3$-degenerate graphs. Their  example has isolation number $1$, 
but it is readily adapted. Start with the subdivision of the complete graph $K_{2k}$ and add one vertex adjacent to all the subdivision vertices.
Then add a perfect matching between the original vertices to yield graph $R_k$. It can be checked that $R_k$ is $3$-degenerate and 
has isolation number $k+1$.


\subsection{Asteroidal-Triple-Free Graphs}

Finally, we consider one more question about the lower packing number. We observed earlier that in most graph classes $\gamma(G) / \rho_L(G)$ is unbounded,
even trees. But we show here that the ratio is bounded for asteroidal-triple-free (AT-free) graphs. 
Bonamy et al.~\cite{BCGY25} showed that for every AT-free graph $G$ it holds that $\gamma(G) \le 3 \rho(G)$.
We show that there is a similar bound involving the lower packing number,
though the constant is not the same.

\begin{theorem} \label{t:asteroid} \ \\
(a) If $G$ is an interval graph, then  $ \gamma(G) \le 3 \gamma_2(G) $, and this is sharp. \\
(b) If $G$ is a permutation graph, then  $ \gamma(G) \le 4 \gamma_2(G) $, and this is sharp. \\
(c) If $G$ is an AT-free graph, then  $ \gamma(G) \le 5 \gamma_2(G) $. 
\end{theorem}
\begin{Proof}
(a) Let $S$ be a distance-$2$ dominating set and 
consider any vertex \mbox{$v\in S$}.
It suffices to show that there is a set of three vertices that dominates
all vertices within distance two of $v$. 
In the interval representation of $G$, there is a neighbor $u_\ell$ of $v$ whose interval stretches leftmost and 
a neighbor $u_r$ of $v$ whose interval stretches rightmost. The set 
$\{u_\ell, v , u_r \}$ thus dominates all vertices within distance~two of $v$. It follows that $\gamma(G) \le 3|S|$.

(b) 
Consider a permutation graph $G$ with vertex set
$\{ v_1, \ldots, v_n \}$ generated by permutation $\phi$.
Let $S$ be a distance-$2$ dominating set and 
consider any vertex $v_i\in S$.
It suffices to show that there is a set of four vertices that dominates
all vertices within distance two of $v_i$. 
Let $a$ be the smallest value such that $\phi(a) \ge \phi(i)$, 
$b$ the largest value such that $\phi(b) \le \phi(i)$,
$c$ the value with $c\ge i$ such that $\phi(c)$ is smallest, and
$d$ the value with $d\le i$ such that $\phi(d)$ is largest.
(Note that not necessarily distinct.)
We claim that the set $X = \{v_a,v_b,v_c,v_d\}$ dominates all vertices within distance~two of $v_i$. 

Each vertex of $X$ is either $v_i$ or a neighbor thereof; 
in particular $v_i$ is dominated by $X$. Let $v_j$ be a neighbor of $v_i$. If $j>i$ then $\phi(j)<\phi(i)$ whence $v_j$ is dominated by $v_a$, 
and if $j<i$ then $\phi(j)>\phi(i)$ whence $v_j$ is dominated by $v_b$. So all neighbors of $v_i$ are dominated by $X$. 
Now consider a vertex $v_k$ at distance two from $v_i$, say with common neighbor $v_j$. 
Without loss of generality $k>i$. Then since $v_k$ is not adjacent to $v_i$ it follows that $\phi(k) > \phi(i)$.
Assume first that $j > i$. Then
$\phi(j) < \phi(i)$, and since $v_k$ is adjacent to $v_j$ but not to~$v_i$, it must be that $k < j$.
By the choice of $b$ it follows that $j\le b$, and thus $v_k$ is adjacent to $v_b$.
Assume second that $j < i$. Then by a similar argument the vertex $v_k$ is adjacent to $v_d$.
It follows that $\gamma(G) \le 4|S|$.

(c) Let $G$ be an AT-free graph. Corneil et al.~\cite{COS97asteroid} showed that such a graph
has a pair of antipodal vertices such that the vertices of any path joining
them is dominating. In particular, $\gamma(G) \le \diam(G)+1$ for such a graph $G$. 
In any graph we have that $\gamma_2(G) \ge \frac{1}{5} (\diam(G)+1)$. (Given any shortest path, a vertex can be 
within distance two of at most five vertices of the path.) Hence $ \gamma(G) \le 5 \gamma_2(G)$.
\end{Proof}

\begin{figure}[h]
\noindent\hfill
\includegraphics{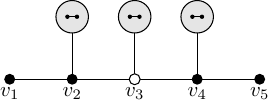} \qquad\qquad
\includegraphics{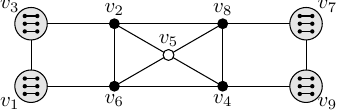}
\hfill\null
\caption{The interval graph $\Itilde$ and permutation graph $\Ptilde$}
\label{f:intervalPermutation}
\end{figure}

The bound of Theorem~\ref{t:asteroid}a is best possible and indeed so is the implication
$\iota(G) \le 3 \rho_L(G)$. 
Construct interval graph $\Itilde$ as follows. Start with a path $v_1 \ldots v_{5}$ on $5$ vertices.
Then for $2 \le i \le 4$ add a  $K_2$ both ends of which are adjacent to $v_{i}$.
The graph $\Itilde$ is drawn on the left in Figure~\ref{f:intervalPermutation}.
Then define $\Itilde_r$ by taking $r$ copies of $\Itilde$ in a row and adding an edge between $v_5$ of one 
copy and $v_1$ of the next.
It is straight-forward to show that $\Itilde_r$ is a connected interval graph with  $\rho_L(I_r)=r$ and $\iota(I_r)=3r$. 

The bound of Theorem~\ref{t:asteroid}b is best possible and indeed so is the implication
$\iota(G) \le 4 \rho_L(G)$. Construct permutation graph $\Ptilde$ as follows.
Start with the permutation graph generated by $(3, 6, 1, 8, 5, 2, 9, 4, 7 )$.
We observe that vertex $v_5$ has eccentricity $2$. 
Further, permutation graphs are closed under cloning vertices both with the new vertex adjacent to the original and with it not.
So, we clone each of the four vertices $v_1$, $v_3$, $v_5$, and $v_7$ to form a copy of $3K_2$.
The graph $\Ptilde$ is drawn on the right in Figure~\ref{f:intervalPermutation}.
(In this example, to dominate all vertices at distance two to vertex $v_5$, we have in our set $X$ that $v_a=v_2$, $v_b=v_8$, $v_c=v_6$, and $v_d=v_4$.)
Then define $\Ptilde_r$ by taking $r$ copies of $\Ptilde$ in a row and 
adjusting the permutation such that the last vertex of one copy is adjacent to the first vertex of the next.
It is straight-forward to show that $\Ptilde_r$ is a connected permutation graph with  $\rho_L(\Ptilde_r)=r$ and $\iota(\Ptilde_r)=4r$. 

It is unclear if the bound of Theorem~\ref{t:asteroid}c is sharp.

\section{Future Work} \label{s:conclude}

For future work, we recall some open questions mentioned earlier:
For a cubic graph $G$ does it hold that $\gamma(G) \le 2 \iota(G)$?
For a connected subcubic graph $G$, except the Petersen graph, does it hold that $\iota(G) \le 2 \rho (G)$?
For a  MOP $G$, does it hold that $\iota(G) \le \rho(G)$?
There are of course multiple other parameters that are related to the five that we focussed on, and 
similar questions can be asked for them.


\end{document}